\documentclass[11pt,a4paper,leqno]{article}
\usepackage{latexsym}
\usepackage{amsfonts}
\usepackage {amsbsy}
\usepackage {amsmath}
\usepackage{amssymb}
\newtheorem{defn}{Definition}[section]
\newtheorem{thm}[defn]{Theorem}
\newtheorem{prop}[defn]{Proposition}
\newtheorem{lem}[defn]{Lemma}
\newtheorem{cor}[defn]{Corollary}
\newtheorem{exa}[defn]{Example}
\newtheorem{remark}[defn]{Remark}
\newcommand \psh{plurisubharmonic }
\newcommand \sx{subextension }
\newcommand \MA{Monge-Amp\`ere }
\newcommand \demo{ Proof: }
\newcommand \C{\mathbb C}
\newcommand \N{\mathbb N}

\newcommand \R{\mathbb R}
\newcommand \B{\mathbb B}
\newcommand \fin{$\blacktriangleright$\\}

\newcommand \lra{\longmapsto}

\newcommand \mc{\mathcal}
\newcommand \mb{\mathbb}
\newcommand \mrm{\mathrm}
\newcommand \Om{\Omega}
\newcommand \om{\omega}
\newcommand \sm{\setminus}

\newcommand \sub{\subset}
\newcommand \Sub{\Subset}

\newcommand \ga{\gamma}

\newcommand \ep{\varepsilon}
\newcommand \f{\varphi}
\newcommand \si{\sigma}

\newcommand \te{\vartheta}
\newcommand \ti {\Tilde}

\newcommand \ove{\overline}
\newcommand \bd{\partial}
\newcommand \we{\wedge}
\numberwithin{equation}{section}
 \title{Maximal subextensions of plurisubharmonic functions}
 \author{U. Cegrell, S. Ko\l odziej, A. Zeriahi}

\begin{document}
\maketitle
\noindent{\bf Abstract}
In this paper we are concerned with the problem of local and global subextensions of (quasi-)plurisubharmonic functions from a "regular" subdomain of a compact K\"ahler manifold.  We prove that a precise bound on the complex Monge-Amp\`ere mass
of the given function implies the existence of a subextension to a bigger regular subdomain or to the whole compact manifold.  In some cases we show that the maximal subextension has a well defined complex Monge-Amp\`ere measure and obtain precise estimates on this measure. Finally we give an example of a plurisubharmonic function with a well defined Monge-Amp\`ere measure and the right bound on its Monge-Amp\`ere mass on the unit ball in $\C^n$ for which the maximal subextension  to the complex projective space $\mb P_n$ does not  have a globally well defined  complex Monge-Amp\`ere measure.

\section {Introduction}
This is the sequel to our earlier paper [CKZ]. There we proved
that given a \psh function $\f $ from the  class $\mc F
(\Om )$ (see the next section for definitions) in a hyperconvex
domain $\Om \Subset \C^n$ one can find its maximal subextension $\ti\f $ which is
\psh in $\C ^n$ and which has logarithmic growth at infinity. If,
in addition, the \MA measure of $\f$ vanishes on pluripolar sets
then the \MA of $\ti\f$ is a well defined positive measure on $\C^n$  in the sense that it is
the weak limit of the sequence of positive measures $(dd^c \f ^j )^n $ for any sequence of continuous
\psh functions $\f ^j \downarrow \ti\f $ having the same rate of
growth at infinity as $\f$. In section 4.3 of this article we
complete this picture studying in more detail the \MA measures of
maximal subextensions $\hat\f$. If the sublevel sets of those
subextensions are bounded then such a measure can be split into
$\mu _1$, dominated by $(dd^c \f )^n $ and essentially supported
on the contact set where $\f =\hat\f$, and $\mu _2$ living on the
set $\partial \{ \hat\f <0 \}.$ In general the maximal global
subextension of  a function from the class $\mc F (\Om )$ may not
have well defined \MA measure. It is the case for generic
multipole Green function as we show in the last section.

Now, a subextension of a \psh function from a domain $D$ in $\C
^n$ to a function defined in  the whole space and of logarithmic
growth can be viewed upon as a subextension of an $\om$-\psh
function (with $\om$ a multiple of the Fubini-Study form) from a
subset of $\C\emph{P} ^n$ to the whole manifold. Here the domain
$D$ is special since there exists a potential for $\om $ in $D.$
If, for instance, $D\subset \C\emph{P} ^n $ contains an algebraic
set of positive dimension then there are no strictly \psh
functions in $D.$ Thus on a compact K\"ahler manifold $X$ we face
a more general problem of \sx of an $\om$-\psh function in $D\subset
X$ to an $\om$-\psh function in $ X.$ In section 3 we
introduce classes of $\om$-\psh functions on $D\subset X$ modelled
on the classes defined by Cegrell and prove the subextension results which
are generalizations to the ones on global subextensions in $\C ^n$.
We refer to [CKZ] for a historical account on subextension
problems.

\section {Monge-Amp\`ere measure of maximal subextensions}

We assume the notational convention $d^c =\frac{i}{2\pi
}(\bar{\partial} -\partial ).$ Let us  recall some definitions
from ([Ce1], [Ce2]). Let $D \Sub \C^n$ be a hyperconvex domain. We
denote by $\mc E_{0} (D)$ the set of negative and bounded \psh
functions $ \f$ on $D$ which tend to zero at the boundary and
satisfy $\int_{D} (dd^c \f )^n < + \infty.$

 Let us denote by $\mc F (D)$ the set of all
$\f \in PSH (D)$ such that there exists a sequence $(\f _{j})$ of
\psh functions in $\mc E_{0} (D)$ such that $\f _{j} \searrow \f$
and $\sup_{j} \int_{D} (dd^c \f _{j})^n < + \infty.$

Before we consider the subextensions from a hyperconvex domain to
$\C ^n$ we first need a result on subextensions to just a larger
hyperconvex set. Let $D\Subset\Omega \Subset \C^n$ be two bounded hyperconvex domains
(open and connected) and and let $u\in \mathcal F(D)$ be a given function. Then $u$ admits  a subextension  $\tilde u \in \mc F (\Omega)$ i.e. $\tilde u \leq u$ on $D$ (see [CZ]). Therefore we can define the maximal subextension of $u$ by 
$$\hat u = \sup\{v\in PSH (\Omega); v<0, \ v|_D \leq u\}.\leqno (\star)$$
It follows from [Ce2] that $\hat u \in \mc F (\Omega)$. 
The following theorem provides a description of the Monge-Amp\`ere
measure of the maximal subextension.

\begin{thm} Let $D\subset\subset\Omega.$
For every  $u\in \mathcal F(D), \hat u\in \mathcal F(\Omega)$,
$(dd^c\hat u)^n \leq \chi_D (dd^cu)^n$  and $\int\limits_{\{\hat
u< u\}}(dd^c\hat u)^n =0.$
\end{thm}

For the proof of the last equality we need the following elementary lemma.
\begin{lem} Suppose $(\mu_j)$ is a sequence of positive measures on $D$
with uniformly bounded mass and that to every $\epsilon > 0$ there
is a $\delta > 0$ such that to every $E \subset D$ with $cap(E) <
\delta$ we have  $\mu_j(E) < \epsilon$  for all j. If $\lim\mu_j =
\mu$ and if $f,g\in PSH(D)$ then

$$\int\limits_{\{f<g\}}d\mu \leq \liminf\limits_j \int\limits_{\{f<g\}}d\mu_j.$$

\end{lem}
To prove the lemma, one can use Bedford-Taylor capacity and the quasicontinuity of $g$ (see \cite{BT2}).

 \noindent\demo(Of the theorem) The first statement of the theorem was proved in [CH].

 Observe that the function $\hat u$ defined by $(\star)$ is plurisubharmonic if $u$ is just any continuous function on $D$. Using the balayage procedure, it is easy to show that in that case we have  $\int\limits_{\{\hat u< u\}}(dd^c\hat u)^n =0$.

   Assume now that $u\in\mathcal F\cap L^{\infty }(D)$ and take a sequence of continuous functions $u_j$ on $D$ decreasing to $u.$ Then $\hat u_j$ decreases to $\hat u$ and the sequence $(\hat u_j)$
 is uniformly bounded on $\Omega$ since $\hat u \leq \hat u_j  \leq 0$ on $\Omega$. Therefore the Monge-Amp\`ere measures
 $(dd^c\hat u_j)^n$ are uniformly dominated by the Monge-Amp\`ere capacity.

 So if we put $\mu_j = (dd^c \hat u_j)^n $ we can apply the lemma
to conclude that for every $s \geq 0$:

 $$
 \int\limits_{\{\hat u_s<u\}}(dd^c\hat u)^n \leq\liminf\limits_j \int\limits_{\{\hat u_s<u\}}(dd^c\hat u_j)^n \leq\liminf\limits_j
 \int\limits_{\{\hat u_j< u_j\}}(dd^c\hat u_j)^n = 0,
 $$
 since by the remark at the beginning of  this proof $\int\limits_{\{\hat u_j< u_j\}}(dd^c\hat u_j)^n = 0.$
 To complete the proof in this case, we let $s$ tend to $+\infty$.

 If $u\in\mathcal F(D)$ only, consider $u_j = \max \{u,-j\}.$ Then, for $t>0$ fixed
 $$
 \left(1 + \max \{u \slash t, -1\}\right) (dd^cu_j)^n \to
 \left(1 + \max \{u \slash t, -1\}\right) (dd^cu)^n, j \to +\infty.
 $$
 Observe that the function $\left(1 + \max \{u \slash t, -1\}\right)$ vanishes on $\{u \leq - t\}$ and is bounded from above by $1$.
 Moreover for any $j > t$ we have $\{u>- t\} \subset \{u > - j\}$ and the sequence of measures ${\bf 1}_{\{u>- j \}}(dd^c u_j)^n$ increases to
 the measure ${\bf 1}_{\{u>-\infty\}} (dd^cu)^n$ (see \cite{BGZ}). Therefore we obtain for $j > t$
 $$
 \left(1 + \max \{u \slash t, -1\}\right) (dd^c u_j)^n \leq {\bf 1}_{\{u>- j\}} (dd^c u_j)^n \leq
 {\bf 1}_{\{u>-\infty\}} (dd^c u)^n.
 $$
 It follows that, for every fixed $t$, the sequence of measures
 $$
 \mu_j := \left(1 + \max \{u \slash t, -1\}\right) (dd^c u_j)^n
 $$
 and therefore $\left(1 + \max \{u \slash t, -1\}\right) (dd^c \hat u_j)^n$ satisfy the requirements of
 the lemma, so we get for every fixed $s$ and $t$:
 $$
 \int\limits_{\{\hat u_s<u\}}\left(1 + \max \{u \slash t, -1\}\right) (dd^c\hat u)^n \leq\liminf\limits_j
 \int\limits_{\{\hat u_s<u\}}\left(1 + \max \{u \slash t, -1\}\right) (dd^c\hat u_j)^n\ $$
 $$
 \leq \liminf\limits_j \int\limits_{\{\hat u_s<u\}}(dd^c\hat u_j)^n\leq\liminf\limits_j
 \int\limits_{\{\hat u_j<u_j\}}(dd^c\hat u_j)^n =0.
 $$
 We now let $t$ tend to $+\infty$. Then since $ 1 + \max \{u \slash t, -1\} \nearrow {\bf 1}_{\{u > - \infty\}}$ as $t \nearrow + \infty$,
 it follows from the previous inequalities that $\int\limits_{\{\hat u_s<u\}}(dd^c\hat u)^n =0.$ To complete the
 proof, we let $s$ tend to $+\infty.$ \fin

\begin{remark} Independently the above theorem was proved in [P], Lemma 4:5.
\end{remark}

\begin{remark} It follows that
$$
 {\bf 1}_{\{\hat u=-\infty\}}(dd^c\hat u)^n = {\bf 1}_{\{u=-\infty\}}(dd^c u)^n.
$$
\end{remark}
Indeed, the inequality "$\leq $" follows from Theorem 2.1 and the
other one from Demailly's inequality [D] (see also [ACCP], Lemma
4.1).

 \section{Potentials on K\"ahler domains}
  Here we want to establish some elementary facts in pluripotential theory on compact K\"ahler manifolds with boundary i.e. on domains in a compact  K\"ahler manifold.

 \subsection{The comparison principle}
 The aim of this section is to give a semi global version of the comparison principle which contains the local one from pluripotential theory on
 bounded hyperconvex domains in $\C^n$ as well as  the global one from the theory on compact
 K\"ahler manifolds (see \cite{GZ2}).

 Let $X$ be a K\"ahler manifold of dimension $n$ and $\om$ K\"ahler form on $X$.
 We want to consider bounded $\om-$plurisubharmonic functions on K\"ahler domains in $X$ with boundary.  For any domain $D \subset  X$, denote by
 $PSH (D,\om)$ the set of $\om-$plurisubharmonic functions on $D$.

 By definition if $\f$ is $\om-$plurisubharmonic on $D$ then locally in $D$ the function $u :=   \f + p$ is a local plurisubharmonic function, where $p$ is a local
 plurisubharmonic potential of the form $\om$ i.e. $dd^c p = \omega$. Therefore the curvature current $\om_{\f} := dd^c \f + \om$ associated to $\f$
 is a globally defined closed positive current on $D$ which can be witten locally as $\om_{\f} = dd^c u$. Therefore by Bedford and
 Taylor [BT], the wedge power $\om_{\f}^p$ is a well defined closed  positive current of bidegree $(p,p)$ on $D$.
  More generally, if $\f_1, \cdots, \f_q$ are bounded $\om-$plurisubharmonic functions on $D$, we can define inductively  the wedge intersection
 product
 \begin{equation}
 T (\f_1, \cdots, \f_q) := \om_{\f_1} \wedge \cdots
 \wedge \om_{\f_q}
 \label{eq:CBT}
 \end{equation}
 as a closed positive current of bidimension $(n- q,n - q)$ on $D$.  Moreover these currents put no mass on pluripolar sets.

  Actually all local results from pluripotential theory
 concerning bounded plurisubharmonic functions on domains in $\C^n$  are valid in the situation considered here. We will
 refer to these results as results from the "local theory".

 Here we use ideas from the global case (see \cite{GZ2}). Our starting point is  the following
 "local version" of the comparison principle which follows from
 quasi-continuity of plurisubharmonic functions (see \cite{BT2},\cite{BT3}).

 \begin{prop}  Let $T$ be a closed positive current of bidimension
 $(p,p)$ ($1 \leq p \leq n$) of type (\ref{eq:CBT})
 and $\f, \psi \in PSH (D,\om) \cap L^{\infty} (D)$. Then
 \begin{equation}
 {\bf 1}_{\{\f < \psi\}} (\om + dd^c \sup\{\f,\psi\})^p \we T = {\bf 1}_{\{\f < \psi\}} (\om + dd^c \psi)^p \we T,
 \label{eq:SCP1}
 \end{equation}
 in the weak sense of Borel measures on $D$.
 In particular
 \begin{equation}
 {\bf 1}_{\{\f \leq \psi\}} (\om + dd^c \sup\{\f,\psi\})^p \we T \geq {\bf 1}_{\{\f \leq \psi\}} (\om + dd^c \psi)^p \we T,
 \label{eq:SCP2}
 \end{equation}
 in the weak sense of Borel measures on $D$.
 \end{prop}

  To perform a useful integration by parts formula, we need to
 consider special domains.
 \begin{defn}
  We will say that a domain $D \subset X$ is quasi-hyperconvex if  $D$
 admits  a continuous negative $\om-$plurisubharmonic exhaustion
 function $\rho : D \lra [- 1 , 0[$.
 \end{defn}
 Observe that any domain $D \subset X$ with smooth boundary given by $D := \{ r <  0\}$,
 where $r$ is smooth in a neighbourhood of $\ove D$, is
 quasi-hyperconvex since
 for $\ep > 0$ small enough, the function $\rho := \ep \ r$ is
 $\om-$plurisubharmonic on a neighbourhood of $\overline D$ and is
 a bounded exhaustion for $D$. Observe that such a domain can be
 pseudoconcave.

Here we will consider only quasi-hyperconvex domains  $D$ satisfying 
 
\begin{equation}
\int_D \om^n <  \int_X \om^n.  
\label{eq:QHC}
 \end{equation}
 \begin{defn} Given a quasi-hyperconvex domain $D$, we define the class of test
 functions $\mc P_0 (D,\om)$ to be the class of functions $\f \in
 PSH^- (D,\om) \cap L^{\infty} (D)$
 such that $\lim_{z \to \bd D} \f = 0$ and $\int_D (\om + dd^c
 \f)^n < + \infty.$
 \end{defn}
 Observe that  for any negative smooth function $h$ with compact
 support in $D$, the function $\ep h$ is in $ \mc P_0 (D,\om)$ for $\ep > 0$ small enough. Moreover, if $\rho$ is an $\om-$plurisubharmonic defining
 function for $D$ then  for any $0 \leq t \leq 1,$ $t \rho \in \mc
 P_0 (D,\om)$.

 \begin{lem} \label{lem:Comp}
 Let $T$ be a closed positive current of bidimension $(p,p)$ ($1 \leq p \leq n$) of type (\ref{eq:CBT}) and $\f, \psi \in PSH (D,\om) \cap
 L^{\infty} (D)$ such that $(\f - \psi)_{\star} \geq 0$ on $\partial D$. Then we have
 $$
 \int_{\{ \f < \psi\}} \om_{\psi}^p \we T \leq  \int_{\{ \f <
 \psi\}} \om_{\f}^p \we T,
 $$
 and
 $$
 \int_{\{ \f \leq \psi\}} \om_{\psi}^p \we T \leq  \int_{\{ \f \leq
 \psi\}} \om_{\f}^p \we T.
 $$
 and  if  $\f \leq \psi$ on $D$ then
 $$
 \int_D \om_{\psi}^p \we T \leq \int_D \om_{\f}^p \we T.
 $$
 In particular if $\f \in PSH^-(D,\om) \cap L^{\infty} (D)$ and $\f \to 0$ at the boundary, then
 $$ \int_D \om^p \we T \leq \int_D \om_{\f}^p \we T.$$
 \end{lem}
 \demo Recall that the condition $(\f - \psi)_{\star} \geq 0$ means  that for any $\ep > 0$, $\{\f < \psi - \ep\} \Subset D$. So
 replacing $\psi$ by $\psi - \ep$ and letting $\ep \searrow 0$, we can assume that $\{\f < \psi \} \Sub D$.
 Then the function $ \te := \sup \{ \f , \psi\} \in PSH (D,\om) \cap L^{\infty} (D)$ coincides with $\f$ near the boundary of
 $D$. This implies that
 \begin{equation}
 \int_{D} (\om + dd^c \te)^p \we T
 = \int_{D} (\om + dd^c \f)^p \we T.
 \label{eq:ST}
 \end{equation}
 Indeed, using local regularization of plurisubharmonic
 functions, we see that $(\om + dd^c
 \te)^p \we T
 - (\om + dd^c \f)^p \we T = d S,$ in the sense of currents on
 $D$, where $S :=   d^c (\te - \f) \left((\om + dd^c \
 \te)^{p-1} + \cdots + (\om + dd^c \f)^{p-1}\right) \we
 T$ is a well defined current with measure
 coefficients and with compact support in $D$. Therefore, by
 definition of the differential of a current, we get $\int_{D} \chi d S = 0$ for any test function $\chi$ which is identically $1$ in q neignbourhood of the support of $S$.  This implies the
 identity (\ref{eq:ST}).

 Now by Proposition 3.1, we get
 $$ \int_{\{ \f < \psi\}} \om_{\psi}^p \we T = \int_{\{ \f < \psi\}} \om_{\te }^p \we T.$$
 Then using  the identity (\ref{eq:ST}) and again Proposition 3.1, we deduce
 \begin{eqnarray*}
 \int_{\{ \f < \psi\}} \om_{\psi}^p \we T & = & \int_D \om_{\te }^p \we T -  \int_{\{ \f \geq \psi\}} \om_{\te }^p \we T  \\
 &  \leq & \int_D \om_{\f}^p \we T  -  \int_{\{ \f > \psi\}} \om_{\te }^p \we T  \\
 & = & \int_D \om_{\f}^p \we T - \int_{\{ \f > \psi\}} \om_{\f }^p \we T,
 \end{eqnarray*}
 which implies
 $$ \int_{\{ \f < \psi\}} \om_{\psi}^p \we T \leq  \int_{\{ \f \leq \psi\}} \om_{\f}^p \we T.$$
 Applying this result to  $\f + \ep$ and $\psi$ and letting $\ep \to 0$, we obtain the required inequality.

 To obtain the second inequality, we can assume $\f, \psi < 0$ on D.
 Now apply the above inequality to $\f$ and $t \psi$ with $0 < t < 1$ and observe that $(dd^c (t \psi) + \om)^n \geq t^n \om^n_{\psi}$. Then letting $t \to 1$, we obtain the required inequality.
 \fin
 If $m = 0$ we set $T_0 = 1$ and for $ m \geq 1$ we set $T_m := \om_{u_1} \we \cdots \we \om_{u_m},$  where $u_1, \cdots u_{m} \in \mc P_0 (D,\om)$. Thus $T_m$ is a closed positive current on $D$.
 Then we have the following important result.

 \begin{cor} 1) The class $ \mc P_0 (D,\om)$ is convex and satisfies the lattice condition:
 $$
 \f  \in \mc P_0 (D,\om), u \in PSH^- (D,\om) \Longrightarrow  \sup \{\f, u\} \in P_0 (D,\om).
 $$
 2) Let $ 1 \leq p,q$ be integers such that $  p + q \leq n$ and denote by $m := n - p - q$. Then for any $\f, \psi  \in P_0 (D,\om),$
 \begin{equation}
 \int_D \om_\f^p \we \om_\psi^{q} \we T_m \leq \int_D \om_{\f}^{p + q} \we T_m +
 \int_D \om_{\psi}^{p + q} \we T_m.
 \label{eq:MMA}
 \end{equation}
 3) If $\f_1, \cdots \f_n \in \mc P_0 (D,\om)$. Then
 $$
 \int_D \om_{\f_1} \wedge \cdots
 \wedge \om_{\f_n}  \leq 2^{n- 1} \sum_{j = 1}^n \int_D \om_{\f_j}^n.
 $$
 \end{cor}
 \demo Let $\f \in \mc P_0 (D,\om)$  and $u \in PSH^- (D,\om)$ and
 denote by $\si (\f,u) := \sup \{\f,u\}$. Since $ \f \leq \si
 (\f,u) \leq 0$, it is clear from the lemma above that
 $$
 \int_D \om^n \leq \int_D (\om + dd^c \si (\f,u))^n \leq \int_D \om_{\f}^n,
 $$
 which implies that  $\si (\f,u) \in P_0 (D,\om)$.

 Now we prove the inequality (\ref{eq:MMA}).

 Indeed, by Lemma 3.4 we get
 $$
 \int_{\{ \f +\epsilon< \psi \}} \om_\f^p \we \om_{\psi}^{q} \we T_m
 \leq \int_{\{ \f < \psi \}} \om_{\f}^{p + q} \we T_m \leq \int_D \om_{\f}^{p + q} \we T_m.
 $$
 Applying this result with $\psi = 0$ we deduce that
 \begin{equation}
 \int_D \om_\f^{p} \we \om^{ q} \we T_m
 \leq \int_{D} \om_{\f}^{p + q} \we T_m.
 \label{eq:MI}
 \end{equation}

 In the same way we obtain
 $$
  \int_{\{  \psi <  \f \}} \om_\f^p \we \om_{\psi}^{q} \we T_m  \leq \int_D \om_{\psi}^{p + q} \we T_m.
 $$
 Therefore
 $$
 \int_{D} \om_\f^p \we \om_{\psi}^{q} \we T_m
 \leq \int_{D} \om_{\f}^{p + q} \we T_m +  \int_{D} \om_{\psi}^{p + q} \we T_m,
 $$
 if we choose $\epsilon > 0$ such that $\int_{\{  \psi + \epsilon = \f \}} \om_\f^p \we \om_{\psi}^{q} \we T_m = 0$
 and let $\epsilon$ decrease to 0.

 The convexity of $\mc P_0 (D,\om)$ follows immediately from the last inequality  since for  $\f, \psi \in \mc P_0 (D,\om)$ and $0 < t < 1$, we have
 $$ \left(\om + dd^c (t \f + (1-t) \psi)\right)^n = \sum_{p = 0}^n {n \choose p} t^p (1 - t)^{n - p} \om_{\f}^p \we \om_{\psi}^{n - p},$$
 which implies by the previous inequality for $m = 0$
 $$
 \int_D \left(\om + dd^c (t \f + (1-t) \psi)\right)^n \leq \int_{D} \om_{\f}^{n} + \int_{D} \om_{\psi}^{n}.
 $$
 To get the last inequality we proceed by induction applying the previous inequality.\fin

 \subsection{Integration by parts formula}

 To prove the integration by parts formula (IBP) which will be crucial for our considerations, we need a semi-global version of the classical (local) convergence theorem of Bedfod and Taylor  for our class $\mc P_0 (D,\om)$.

 \begin{prop} Let
 $(\f_j^0), \cdots (\f_j^n)$ be sequences of locally  uniformly
 bounded $\om-$plurisubharmonic functions in the class $\mc P_0 (D,\om)$  converging monotonically to  $\f^0, \cdots, \f^n \in \mc P_0 (D,\om)$ respectively.
 Then the  positive currents $S_j := (dd^c \f_j^1 + \om) \wedge \cdots \wedge (dd^c
 \f_j^n + \om)$ and $S := (dd^c \f^1 + \om) \wedge \cdots \wedge (dd^c
 \f^n + \om)$ have uniformly bounded total masses in $D$ and
 $$
 \lim_{j \to + \infty} \int_D (-\f_j^0) (dd^c \f_j^1 + \om) \wedge \cdots \wedge (dd^c
 \f_j^n + \om)  =
 $$
 $$
 \int_D (-\f^0) (dd^c \f^1 + \om)
 \wedge \cdots \wedge (dd^c \f^n + \om).
 $$
 \end{prop}
 \demo
 Observe first that the local theory of Bedford and Taylor implies that $(-\f_j^0) S_j \to
 (-\f^0) S$ weakly on $D$ (see [BT2]).  It
 follows from our hypothesis that given $\ep  > 0$, there exists
 an open set $D' \Sub D$ such
 that  $- \ep \leq \f_j^0 \leq 0$ and $- \ep \leq  \f^0 \leq 0$ on
 $D \sm D'$.  Then
 \begin{equation}
 \int_{D} (-\f_j^0) S_j - \int_{D} (-\f^0) S  =  \int_{D'}
 (-\f_j^0) S_j - \int_{D'} (-\f^0)
 S_0 \ + \ O (\ep),
 \label{eq:SE}
 \end{equation}
 uniformly in $j \in \N$. Here we have used the fact that the
 currents $S_j$ have uniformly bounded mass on $D$ by Lemma \ref{lem:Comp}.
 Now observe that we can always choose the domain $D'$ so that the
 positive measure $\mu_0 := (- \f^0) S$
 puts no mass on its boundary $\bd D'$. Then since the positive
 measures $\mu_j := (-\f_j^0) S_j $ converge
 weakly to $\mu_0$ in $D$, it follows that
 $$
 \mu_0 (D') \leq
 \liminf_j \mu_j (D') \leq  \limsup_j \mu_j (\ove{D'}) \leq \mu_0
 (\ove{D'}) = \mu_0 (D'),
 $$
 which proves that the first integral on the right hand side
 converges to $0$ and the proposition is proved.

 Now we can prove the following integration by parts formula which
 will be  useful in the sequel.
 \begin{lem} Let $T := (\om + dd^c u_1) \we \cdots \we (\om + dd^c u_{n - 1})$, where $u_1, \cdots u_{n - 1} \in \mc P_0 (D,\om)$. Let
  $u, v\in \mc P_0 (D,\om)$. Then
 \begin{equation}
  \int_D u dd^c v \we T = \int_D v dd^c u \we T,
 \label{eq:IBP1}
 \end{equation}
 and
 \begin{equation}
  \int_D u \om_v \we T - \int_D v \om_u \we T = \int_D (u - v) \om \we T.
 \label{eq:IBP2}
 \end{equation}
 \end{lem}
  \demo
 Denote by $H (u,v) := u dd^c v \we
 T -  v dd^c u \we T$.  Then by Proposition 3.1 the current $H (u,v)$  has finite total mass in $D$. It follows from Stokes formula that
 if $\bar u, \bar v$ are  bounded $\om-$plurisubharmonic functions
 on $D$ such that $\bar u = u$ and
 $ \bar v = v $ near the boundary $\bd D$ then
 $$\int_D H (\bar u,\bar v) = \int_D H (u,v).$$
 Indeed observe that since $u, v, \bar u, \bar v$ are bounded
 $\om-$quasiplurisubharmonic functions on $D$, it follows from the
 local theory that the currents $S  := u d^c v \we T  - v d^c u
 \wedge T$ and $\bar S := := \bar u d^c \bar v \we T  - \bar v d^c
 \bar u \wedge T$ are
 well defined currents with measure coefficients on $D$ such that
 $ d S =  u dd^c  v \we
 T -  v dd^c u \we T = H (u,v)$ and $d \bar S = \bar u dd^c \bar v
 \we T  - \bar v dd^c \bar u \wedge T = H (\bar u,\bar v)$ in the
 weak sense of currents on $D$. Now since $S - \bar S$ is of
 compact support in $D$, it follows that $\int_D d (S - \bar S) =
 0$  and then
 $$\int_D H (u,v) = \int_D H (\bar u,\bar v).$$
 Now for $\ep > 0$ small enough, set
 $u_{\ep} := \sup \{u , v - \ep\}$  and $v_{\ep} := \sup \{v , u -
 \ep\}$ and observe that $u_{\ep} = u$ and $v_{\ep} = v$ near $\bd
 D$. Thus by the previous remark, we have for $\ep > 0$
 small enough
 \begin{equation}
 \int_D H (u_{\ep},v_{\ep})  =  \int_D H (u,v).
 \end{equation}
 We want to pass to the limit. Here we must use the fact that $u =
 v = 0$ on $\bd D$, which implies that $u_{\ep} = v_{\ep} = 0$ on
 $\bd D$.
 Now for $\ep > 0$ small enough, we have
 $$ H (u_{\ep},v_{\ep}) = u_{\ep} dd^c v_{\ep} \we T -  v_{\ep}
 dd^c u_{\ep} \we T.$$
 Since $u_{\ep} \nearrow g := \max\{u , v\}$ and $v_{\ep}
 \nearrow g$, it follows from Proposition 3.1 that
 $$
  \lim_{\ep \to 0} \int_D  u_{\ep} dd^c v_{\ep} \we T =
 \int_D g dd^c g \we T = \lim_{\ep \to 0} \int_D  v_{\ep}
 dd^c u_{\ep} \we T,
 $$
 which implies the required integration by parts formula.
 \fin

 \section{Subextension of quasi-plurisubharmonic functions} 

\subsection{Weighted Monge-Amp\`ere energy classes}
In the contrast to the local case, the domain of definition of the complex Monge-Amp\`ere operator is not well understood in the global case. Interesting classes have been investigated in \cite{GZ2} and \cite{CGZ}.
We are going to introduce similar classes in the semi-global case where the complex Monge-Amp\`ere operator is
well defined and continuous under deacreasing sequences.
The first class is modeled on the class defined by Cegrell in (\cite{Ce2}) as  follows.
 \begin{defn}
 We say that $\f \in \mc F (D,\om)$ if there exists a
 decreasing sequence $(\f_j)$  from the class
 $\mc P_0 (D,\om)$ which converges to $\f$ on $D$ such that
 $$
  \sup_{j}  \int_D \om_{\f_j}^n < + \infty.
 $$
 \end{defn}

 Observe that $\mc F (D,\om)$ is a convex set and $\mc P_0 (D,\om)
 \sub\mc F (D,\om)$. The class  $\mc F (D,\om)$ is the counterpart of the class defined by Cegrell in  \cite{Ce2}. Let
 $D$ be a hyperconvex domain where the form $\om$ has a
 plurisubharmonic potential $q$ on $D$ with boundary values $0$ and let $\mc F (D)$ the  class defined in \cite{Ce2}. Then if $\f \in \mc F (D,\om)$ iff
$u := \f + q
 \in \mc F (D)$.

We do not know at the moment if the Monge-Amp\`ere operator is well defined on the class 
$\mc F (D,\om)$ but we can define the Monge-Amp\`ere mass of a function $\f \in \mc F (D,\om)$ thanks to the following lemma.
\begin{lem} Let $\f \in \mc F (D,\om)$ be a fixed function. Then the constant
$$M_D (\f) := \lim_j \int_D (\om + dd^c \f_j)^n = \sup_j \int_D (\om + dd^c \f_j)^n$$
is independant of the decreasing sequence $(\f_j)$ from $\mc P_0 (D,\om)$ converging to $\f$.

Moreover if    $\psi \in PSH (D,\om)$ and $\f \leq \psi \leq 0$ then $\psi \in \mc F (D,\om)$.

\end{lem} 
\demo Take a defining sequence $(\f_j)_j$ for $\f$.
By Lemma 3.4 we know that the sequence $\{\int_D (\om + dd^c \f_j)^n\}_j$ is increasing and by definition it is bounded so the limit $M_D (\f)$ exists. We only need to show that it does not depend on the sequence. Let $(\psi_j)$ another decreasing sequence of functions in the class $\mc P_0 (D,\om)$ converging to $\f$ in $D$.
 Fix $\varepsilon >0$ and $j$. Since by Bedford-Taylor continuity theorem ([BT2]), $ (\om + dd^c \sup \{\psi_j,\f_{k}\})^n \to (\om + dd^c \psi_j)^n$ weakly on $D$ as $k \to \infty$, it follows  that there exists $k_j$ such that
$$
\int_D (\om + dd^c \sup \{\psi_j,\f_{k_j}\})^n > \int_D (\om + dd^c \psi_j)^n - \varepsilon.
$$
By  Lemma 3.4, we have 
$$
\int_D (\om + dd^c \sup \{\psi_j,\f_{k_j}\})^n \leq \int_D (\om + dd^c \f_{k_j})^n
\leq M_D (\f).
$$
Therefore it follows that $ \int_D (\om + dd^c \psi_j)^n - \varepsilon \leq  M_D (\f)$,
which implies that $ \sup_j \int_D (\om + dd^c \psi_j)^n \leq  M_D (\f)$ and proves the first part of the lemma.

Now set $\psi_j := sup\{\psi, \f_j\}$. Then by  Lemma 3.4, $\psi_j \in \mc P_0 (D)$ and $ \int_D (\om + dd^c \psi_j)^n \leq \int_D (\om + dd^c \f_{j})^n \leq M_D(\f)$.
Since $(\psi_j)$ decreases to $\psi$, it follows that $\psi \in \mc F (D,\om)$ and from the first part of the proof we deduce that $M_D(\psi)  \leq M_D(\f)$.

\fin

 Let us introduce  the following classes of finite weighted Monge-Amp\`ere
 energy (see \cite{Ce1}, \cite{GZ2}, \cite{BGZ}).
 A weight function is by definition an increasing function $\chi : \R \lra \R$ such  that $\chi (t) = t$ is $ t \geq 0$ and  $\chi (- \infty) = - \infty $. To any weight function we associate
 the  class $ \mc E_{\chi} (D,\om) $ of of $\omega-$plurisubharmonic functions $ \f \in PSH (D,\om)$ for which there exists a sequence $(\f_j) \in \mc P_0 (D,\om)$, $\f_j \searrow \f$ such that
 $$ \sup_j \int_D
 \vert \chi (\f_j)\vert \om_{\f_j}^n < + \infty.
 $$
 In our case the weight function $\chi$ will be convex.
From the  (IBP) formula, we can derive the following fundamental inequality which will be useful (see \cite{GZ2}).
 \begin{prop} Let $\chi : \R \lra \R$ be a convex weight function. Then for any $\f, \psi \in \mc P_0 (D,\om)$ with $\f \leq \psi$, we have
 \begin{equation}
 \int_D \vert \chi (\psi)\vert \om_{\psi}^n \leq 2^n \int_D \vert \chi (\f) \vert \om_{\f}^n.
 \label{eq:F-INEQ}
 \end{equation}
 \end{prop}
We can prove that the complex Monge-Amp\`ere operator is well defined and continuous on decreasing sequences in the class $ \mc E_{\chi} (D,\om),$
where $\chi$ is a convex increasing functions $\R\lra \R$ (see \cite{GZ2}, \cite{CGZ}).

 \begin{prop}  The complex  Monge-Amp\`ere operator is well defined on the class $\mc E_{\chi} (D,\om)$. Moreover if $(\f_j)$ is a decreasing sequence from the class
 $\mc E_{\chi} (D,\om)$ which converges to $\f \in \mc E_{\chi} (D,\om)$, then  the Monge-Amp\`ere measures $(\om^n_{\f_j})$
 converge to $\om_{\f}^n$ weakly on $D$. Moreover for any $h\in PSH(D, \omega ) \cap L^{\infty} (D)$
 $$\lim _j \int _D h\om _{\varphi _j} ^n  =\int _D h\om _{\varphi } ^n .$$
 \end{prop}
 
 Using the integration by parts formula, the fundamental inequality and following the same
 arguments as \cite{GZ2}, it is possible to prove the
 following result.
 \begin{prop} Let $\f \in \mc PSH
 (D,\om)$. Assume there exists  a decreasing sequence  $(\f)_{j \in \N}$ in $\mc P_0 (D,\om)$
 which converges to $\f \in PSH (D,\om)$ and satisfies  $\sup_j  \int_D \vert
 \chi (\f_j)\vert \om_{\f_j}^n < + \infty$. Then $\f \in \mc E_{\chi}
 (D,\om)$ and
 $$ \lim_{j \to + \infty} \int_D \vert
 \chi (\f_j)\vert \om_{\f_j}^n = \int_D \vert
 \chi (\f)\vert \om_{\f}^n.$$
 \end{prop}

 \subsection{A general subextension theorem}
 We now prove the following general subextension result which generalizes our previous result with a new proof (see \cite{CKZ}).

 \begin{thm} \label{thm:SUBEXT}
 Let $D \subset X$ be a quasi-hyperconvex domain satisfying the condition $(3.4)$. Let
 $\f \in \mc F (D,\om)$ such that $M_D (\f) \leq \int_X \om^n$.
 Then there exists a function $\ti {\f} \in
 PSH (X,\om)$ such that $\ti \f \leq \f$ on $D$.
 \end{thm}
 \demo   Let $(\f_j)$ be a decreasing sequence from the class
 $\mc P_0 (D,\om)$ which converges to $\f$ on $D$. By Lemma 4.2 we have
 $$
 \int_D (\om + dd^c \f_j)^n \leq M_D (\f).
 $$
 
 First assume that $M_D (\f) < \int_X \om^n$.
 Then by [GZ2]
 there exists $u_j \in \mc E^1 (X,\om)$ with
 $\sup_X u_j  = - 1$ such that
 $$
 (\om + dd^c u_j)^n =  {\bf 1}_D (\om + dd^c \f_j)^n + \ep_j \om^n
  $$
 on $X$, where $\ep_j > 0$ is chosen so that the total mass of both  sides are equal.
 Fix $j \in \N$. Since $ \{\f_j < u_j\} := \{ x \in D ; \f_j < u_j \} \Subset D,$ and $\f_j$ is bounded, it follows that for $s > 1$ large enough,
 $\{\f_j < u_j^s\} = \{\f_j < u_j\} \Subset D$, where $u_j^s := \sup \{u_j, - s\}$. Then
 by the comparison principle (Lemma 3.4), it follows that
 $$
  \int_{\{\f_j < u_j^s\}} (\om + dd^c u_j^s)^n \leq \int_{\{\f_j <
 u_j^s\}} (\om + dd^c \f_j)^n.
 $$
 Recall that ${\bf 1}_{\{ u_j > - s\}} (\om + dd^c u_j^s)^n = {\bf 1}_{\{ u_j > - s\}} (\om + dd^c u_j)^n$ (see \cite{GZ2}). Therefore
 $$
 \int_{\{\f_j < u_j\}} (\om + dd^c u_j)^n \leq \int_{\{\f_j <
 u_j\}} (\om + dd^c \f_j)^n,
 $$
 which implies that $Vol_{\om} (\{\f_j < u_j\}) = 0$ and then $u_j
 \leq \f_j$ on $D$. Due to the normalization of $u_j$, the function  $u := (\limsup_{j \to + \infty} u_j)^* \in PSH (X, \om)$ and
 satisfies $u \leq \f$ on $D$.

 Now assume $\f \in \mc F (D,\om)$ with $M_D (\f) = \int_X \om^n$ and consider a decreasing sequence $(\f_j)$ in  $\mc P_0 (D,\om)$ converging to $\f$ with uniformly bounded Monge-Amp\`ere masses.
 Then  it follows that for any $0 < t < 1$
 the function $t \f_j \in \mc P_0 (D,\om)$ and $\int_D (\om + dd^c t
 \f_j)^n = \int_D (t \om_{\f_j} + (1 - t) \om)^n $.
By  Lemma 3.4 we have $ \int_D \om_{\f_j}^p \wedge\om^{n-p} \leq  \int_D \om_{\f_j}^n$.  Therefore since $\int_D \om^n < \int_X \om^n$, it follows that $M_D (t \f_j) = \int_D (\om + dd^c t
 \f_j)^n  < \int_X \om^n$.
 By the first part we can find a subextension $\psi_{j}^t \in PSH
 (X,\om)$ of $t \f_j$ satisying $\max_X \psi_{j}^t = - 1$. Therefore the
 function $\psi_j  := (\limsup_{t \nearrow 1} {\psi_j}^t)^*$ is an
 $\om-$plurisubharmonic subextension of $\f_j$ to $X$ with $\max_X \psi_{j} = - 1$ . Now observe that $(\psi_j)$ is a decreasing sequence of plurisubharmonic functions on $X$ which converges to  a plurisubharmonic function $\psi$ on $X$ such that $\max_X \psi = - 1$ and $\psi \leq \f$ on $D$.
 \fin

  It follows from the above theorem that given $\f \in \mc F
 (D,\om)$ such that $M_D (\f) \leq \int_X \om^n$, the following function
 $$
 \ti \f = \ti \f_D := \sup \{\psi  \in PSH (X,\om) ; \psi \leq \f \
 \mathrm{on} \ D \}
 $$
 is a well defined $\om-$plurisubharmonic function on $X$ and will be called the maximal subextension of $\f$ from $D$ to $X$.
 
  The example below shows that in general the maximal subextension
 does not belong to the global domain of definition of the complex
 Monge-Amp\`ere operator on $X$ since it may have positive Lelong number along a hypersurface.

 However if the given function has a finite weighted Monge-Amp\`ere  energy in the sense of [GZ2], we will prove that the maximal subextension satisfies the same property.

 \begin{thm} Let $D \subset X$ be an quasi-hyperconvex domain satisfying the condition $(3.4)$ and let $\f \in \mc E_{\chi} (D,\om)$ be such that $\int_D \om_\f^n \leq \int_X \om^n$, where
 $\chi : \R \lra \R$ is a convex weight function. Then the maximal subextension $\tilde {\f}$ of $\f$ from $D$ to $X$ exists and has the following properties:\\
 $(i)$  $\ti {\f} \in \mc E_{\chi}
 (X,\om)$ and $\int_X \vert \chi\circ  \ti {\f} \vert (\om + dd^c \ti
 {\f})^n \leq \int_D \vert \chi\circ {\f}\vert (\om + dd^c {\f})^n$, \\
 $(ii)$  ${\bf 1}_D (\om + dd^c \ti \f)^n \leq {\bf 1}_D (\om +
 dd^c \f)^n$ holds
 in the sense of measures on $X$, \\
 $(iii)$ the measure $(\om + dd^c \ti \f)^n$ is carried by the Borel set $\{\ti
 \f = \f \} \cup \bd D$.
 \end{thm}

  We will need the following lemma which can be proved using the argument from the first part of the proof of Theorem 2.1.
 \begin{lem}
 Let $D$ be as above and $\f \in \mc P_0 (D,\om)$ be such that $\int_D \om_\f^n \leq \int_X \om^n$, then $\ti{\f} \in PSH (X,\om) \cap L^{\infty} (X)$ and
 ${\bf 1}_D (\om + dd^c \ti{\f})^n \leq {\bf 1}_D (\om + dd^c
 \f)^n$ in the sense of measures on $X$. Moreover the measure $(\om
 + dd^c \ti{\f})^n$
 is  carried by the Borel set  $ \{x \in \bar{D} ; \tilde {\f} (x) = \f (x)\}$.
 \end{lem}

 \noindent Proof of the theorem.  Let $(\f_j)$ a sequence $(\f_j) \in \mc P_0 (D,\om)$ which
 decreases to $\f$ on $D$. Define $\tilde {\f}_j$ to be the maximal
 subextension of $\f_j$ from $D$ to $X$. Then by the previous lemma $\tilde {\f}_j
 \in PSH (X,\om) \cap L^{\infty} (X)$  and $(\om + dd^c
 \ti{\f}_j)^n$ is supported
 on the contact set $\{x \in \bar{D} : \tilde {\f}_j (x) = \f_j (x)\}$.
 Hence $(- \chi \circ  \tilde {\f}_j) (\om + dd^c \tilde{\f}_j)^n \leq
 {\bf 1}_D (- \chi \circ \f_j) (\om + dd^c \f_j)^n$ in the sense of
 measures on $X$. Therefore there is a uniform constant $ C> 0$
 such that for any $j \in \N$,
 $$
 \int_X (- \chi \circ \tilde {\f}_j) (\om + dd^c \tilde{\f}_j)^n
 \leq \int_D (- \chi \circ \f_j) (\om + dd^c \f_j)^n \leq C.
 $$
 Since $(\tilde {\f}_j) \searrow \ti \f$ on $X$ it follows from \cite{GZ2} that $\ti \f \in \mc E_{\chi} (X,\om)$. Moreover by the
 convergence theorem (\cite{GZ2}, \cite{CGZ}) it follows that
 ${\bf 1}_D \vert \chi \circ \tilde {\f}\vert (\om + dd^c \tilde{\f})^n
 \leq {\bf 1}_D \vert \chi \circ \f\vert (\om + dd^c \f)^n$ in the
 sense of measures on $X$.

 The third part of the theorem is proved along the same lines as the last part of the proof of Theorem 2.1 using Lemma 4.8 and Proposition 4.2.
 \fin
 
\begin{remark} In contrast to the local case it may happen that a part of the  Monge-Amp\`ere measure of $\tilde \f $ lives on the boundary of $D$. 
\end{remark}
 As we already said before, the example in the last section shows that the maximal subextension of a given function $\f \in \mc F (D,\om)$ may have not a well defined Monge-Amp\`ere measure. However the following property may be useful.

 \begin{prop} Let $\f \in \mc F (D,\om)$ be a given function. Then  if $(\f_j)$ is a decreasing sequence of functions in
 the class  $\mc P_0 (D,\om)$ converging to $\f$ then the sequence $(\ti \f_j)$
 decreases to $\ti \f$ on $X$.
Moreover any Borel measure $\mu$ on $X$ which is a limit point of the sequence of measures ${(\om + dd^c \ti \f_j)^n}$ on $X$ satisfies the inequality ${\bf 1}_D \mu \leq {\bf 1}_D (\om + dd^c \f)^n$ in the sense of measures on $X$.
 \end{prop}
 \demo Observe that for each $j \in \N$, $\ti \f$ is a global
 subextension of $\f_j$ to $X$ and then $ \ti \f \leq \ti \f_j$ on
 $X$. Therefore it is clear that the sequence $(\ti \f_j)$
 decreases to an $\om-$plurisubharmonic function $\psi$ on $X$
 which satisties the inequality $\ti \f \leq \psi$ on $X$. This
 shows that $\psi \in \mc PSH (X,\om)$. On the other hand since
 $\psi \leq \ti \f_j \leq \f_j$ on $D$ we infer that $\psi \leq \f$  on $D$, which proves that $\psi$ is a subextension of $\f$ to $X$
 and then $\psi \leq \ti \f$ on $D$. We conclude that $\psi = \ti
 \f$ on $X$. 
We know from the last lemma that ${\bf 1}_D (\om + dd^c \ti \f_j)^n  \leq {\bf 1}_D (\om + dd^c \f_j)^n$ in the sense of measures on $X$, which implies the last statement of the proposition.
\fin

\subsection{Subextension in $\mathbb C^n$}

 \bigskip
 Now we pass to  subextensions from a hyperconvex domain $D \Subset \mathbb C^n$ to
 $\mathbb C^n,$ considered as an open subset of $\mb P_n$. Recall that the Lelong class is defined by
$$
\mathcal L(\C^n) := \{u\in PSH(\mathbb C^n);
 \sup \{u(z) - \log^+|z| < + \infty  \}.
$$ 
 Let $\omega = \om_{FS}$ be the normalized Fubini-Study metric on $\mb P_n$ defiend in affine ccordinates by
 $$
  \omega := dd^c \log \vert \zeta\vert,
 $$
 where $\zeta := [\zeta_0: \cdots : \zeta_n]$ is the homogenuous coordinates on $\mb P$.
 As usual we will consider $\C^n = \mb P \sm \{\zeta_0 = 0\}$ whith the affine coordinates defined as by $z_j := \zeta_j \slash \zeta_0$ ($1 \leq j  \leq n$).
 With these notations we have $\omega|\C^n = dd^c \ell,$ where $\ell (z) := (1/2) \log (1+|z|^2)$. Therefore given any $u \in \mathcal L (\C^n)$, the 
 function defined by
$$
\f (\zeta) := u (z) - (1/2) \log (1+|z|^2), \zeta_0 \neq 0
$$
is $\omega-$plurisubharmonic on $\mb P_n \sm \{\zeta_0 = 0\}$ and locally upper bounded in a neighbourhood of the hyperplane at infinity $H_{\infty} := \{\zeta_0 = 0\}$ so that it extends to an $\omega-$plurisubharmonic function on $\mb P_n$ which we also denote by $\f$. It follows that the correspondance $u \longmapsto \f$ is a bijection between
$\mathcal L(\C^n)$ and $PSH (\mb P,\om)$ such that $\omega + dd^c \f = dd^c u$ on $\C^n$.

 From the last theorem we can deduce a generalization of our earlier result (see \cite{CKZ}, Theorem 5.3).
 \begin{thm} Let $D \Subset \C^n$ be a hyperconvex domain and let $u \in \mc F (D)$ be such that $(dd^c u)^n$ does not put any mass on pluripolar sets in $D$ and $\int_D (dd^c u)^n \leq 1$. Then its maximal subextension $\tilde   u$ from $D$ to $\C^n$  belongs to $ \mc L (\C^n)$ and has a well defined global Monge-Amp\`ere measure $(dd^c \tilde u)^n $ which is carried by the set
 $\{\tilde  u = u\} \cup \partial D$ and satsifies the inequality  $  {\bf 1}_{D} (dd^c \tilde u)^n \leq {\bf 1}_{D} (dd^c u)^n$.
 \end{thm}
 \demo Assume first that $D = B_R$ is an euclidean ball with center at the origin and radius $R > 0$. Then the function  $q := (1/2) \log (1+|z|^2) - (1/2) log (1+R^2)$ is a potential of the normalized Fubini-Study form $\om$ on $\C^n$ which vanishes on  $\partial D$. In this case $\f := u - q \in \mc F (D,\om)$.  From our hypothesis  $(\om + dd^c \f)^n (\{\f = - \infty\}) = (dd^c u)^n (\{u = - \infty\}) = 0$. It follows from standard fact in measure theory that there
 exists a convex inceasing function $\chi : ]- \infty , 0] \longrightarrow ]- \infty , 0]$ such that $ \int_D (-\chi \circ \f) (\om + dd^c \f)^n < +
 \infty$ (see [GZ 2]. It easily follows that $\f \in \mc E_{\chi} (D,\omega)$ and  then we can apply the last result to find a subextension $\tilde \f \in \mc E (\mb P_n,\omega)$ of $\f$ to $\mb P_n$. Then $\tilde u :=
\tilde \f + q$ is the maximal subextension of $u$  to $\C^n$. 

Now in the general case consider an euclidean ball $B$ such that $D \subset B$ and use Theorem 2.1 to produce a subextension $v \in \mc F (B)$ of $u$. Then by the previous case $v$ has a subextension $\tilde v$ such that $\psi :=\tilde v - q$ is a function in $\mc E (\mb P_n,\omega)$ which is a subextension of $\f := u - q$ from $D$ to $\mb P_n$.
Therefore the maximal subextension $\tilde \f$ of $\f$ exists and since $\psi \leq \tilde \f $ it follows that $\tilde \f \in \mc E (\mb P_n,\omega)$. Thus $\tilde u := \ti \f + q \in \mc L (C^n)$ is the maximal subextension of $u$ to $\C^n$. The other properties follow in the same way as in the proof of Theorem 4.8.
 \fin
  
Now we consider an arbitrary function  $u\in \mathcal
 F(D)$ and a positive $\gamma$  satisfying 
$$
\gamma^n \geq   \int\limits_D (dd^cu)^n.
$$ 
Then from Theorem 4.6 the set of entire subextensions of
 logarithmic growth
 $$
 \{v\in PSH(\mathbb C^n); v|_{D} \leq u, v(z)\leq a_v + \gamma log^+|z| \}$$
is not empty. Thus, using notation
$$
\mathcal L_{\gamma} (\C^n) = \{v\in PSH(\mathbb C^n);
v(z)\leq a_v + \gamma log^+|z| \}
$$
one can choose the maximal subextension of $u$ of  logarithmic growth related to $\gamma$
$$
\hat{u}_{\gamma} =  \sup \{v\in \mathcal L_{\gamma} (\C^n) ;v|_{D} \leq u \}.
$$
 As we shall see the \MA measure of this subextension may not
 exist. If it exists however, one can deduce some information on
 the support of such measure.
 
 Define 
 $$N_u = \{z\in\mathbb C^n; \hat{u}_{\gamma} < 0\}.$$

 \begin{prop}
 Assume that $u\in\mathcal F(D)$ and  let $\gamma^n = \int\limits_D (dd^cu)^n.$ Then for any sequence
 $u_j\in \mathcal E_0(D)\cap C(\bar D),$ decreasing to $u$ if $\mu$ is an accumulation point of $(dd^c \hat{u}_{j,\gamma})^n$ then $\mu = f
 (dd^cu)^n + \nu$ where $0\leq f \leq 1$ is a function vanishing outside $D$ and where $\nu$ is a positive measure, $\rm supp\, \nu
 \subset\partial N_u.$  
 \end{prop}
\demo Assume first that $u\in \mathcal E_0(D)\cap C(\bar D).$ Then $\hat{u}_{\gamma}$ is continuous and  the zero sublevel
 set of $\hat{u}_{\ga }, N_u$ is hyperconvex.

 By definition,  $D\subset N_u$ and by Theorem 5.1 in [CKZ] $D$ is
 not relatively compact in $N_u$. There are two possibilities:

 1) $ D = N_u.$

 2)  $D \neq N_u \subset\subset \mathbb C^n.$

 If 1) occurs then $\hat{u}_{\gamma}$ extends $u$ to a function in $\mathcal L_{\gamma}\cap L^{\infty}_{loc}$ and 
 $${\bf 1}_{N_u}(dd^c
 \hat{u}_{\gamma})^n = {\bf 1}_D (dd^c \hat{u}_{\gamma})^n={\bf 1}_D (dd^cu)^n.
 $$ In particular, if $\gamma^n = \int\limits_D (dd^cu)^n$
 then $(dd^c \hat{u}_{\gamma})^n ={\bf 1}_D (dd^cu)^n$ on $\mathbb C^n.$

 Generically we have 2).  Then on $N_u$, $\hat{u}_{\gamma} $ is equal to $ \hat u$, the maximal local subextension of $u$ from $D$
 to $N_u$. Consider $D_j\subset\subset D_{j+1}\subset\subset D$ an exhaustion sequence  of $D.$ Denote by $\hat u_j$ the
 corresponding local maximal subextension to $N_u$ of the solution
 $u_j\in\mathcal E_0(D)$ to $(dd^c u_j)^n = {\bf 1}_{D_{j-1}}(dd^cu)^n.$ Then $\hat u \leq \hat u_j$ and
 $(dd^c\tilde u_j)^n \leq {\bf 1}_{D_{j-1}}(dd^cu)^n$ on $N_u$ by Theorem 2.1 and so $(dd^c\hat u)^n \leq {\bf 1}_D (dd^cu)^n$ on
 $N_u.$

 Therefore, $(dd^c \hat{u}_{\gamma})^n = f (dd^cu)^n + \nu$ where $0\leq f \leq 1$ is a function vanishing outside $D$ and where
 $\nu$ is a positive measure, $\rm supp\, \nu \subset\partial N_u \cap\partial D.$

 Now consider the general case. Choose a deacreasing sequence $(u_j)$ in $\mathcal E_0(D)\cap C(\bar D),$ decreasing to $u.$ Then $\hat{u}_{j,\gamma}$
 decreases to $\hat{u}_{\gamma}$ and $(dd^c \hat{u}_{j,\gamma})^n = f_j (dd^c u_j)^n + \nu_j$ where $0\leq f_j \leq 1$ is a function
 vanishing outside $D$ and where $\nu_j$ is a positive measure, $\rm supp\, \nu_j \subset\partial N_{u_j}.$ Also $\int (dd^c
 \hat{u}_{j,\gamma})^n = \gamma ^n.$ So if $\mu$ is any weak limit of $(dd^c \hat{u}_{j,\gamma})^n,$ then $\mu = f (dd^cu)^n + \nu$
 where $0\leq f \leq 1$ is a function vanishing outside $D$ and
 where $\nu$ is a positive measure carried by $\partial N_u.$  \fin

 \begin{cor} If, for $u\in\mathcal F(D),$ the set $N_u$ is bounded then the \MA measure of $u_{\gamma}$ is well defined and equal to
 the limit of $(dd^c \hat{u}_{j,\gamma})^n.$
 \end{cor}

 If  $N_u$ is not a bounded hyperconvex set, $u_{\gamma}$ need not to be in the domain of definition of the \MA operator.
 This is shown in the following example.

 \begin{exa} The maximal entire subextension of a
 function from the class $\mc F (\B )$ may not have well defined global Monge-Amp\`ere measure on $\C^2 $.
 \end{exa}

 Consider the Green function $g$ in the ball $\B (0,2)\subset \C ^2$ with two poles at $(-1,0)$ and $(1,0)$ of weight
 $\frac{1}{\sqrt{2} }$ each. Then
 $$
 \int _{\B (0,2)} (dd^c g)^2 =1.
 $$
 So there exists the maximal entire subextension $\hat g = \hat g_t$ in the Lelong class $\mc L _t(\C ^2 ),1\leq t < \sqrt 2.$ Note that  $\frac{1}{\sqrt{2}
 }\log ||\frac{ z_2}{2} ||$ is a subextension.
 By the definition of the Green
 function we have for some $R\in (0,1), \ A>0$ the following
 inequalities
 \begin{eqnarray}
 &|g(z) - \frac{1}{\sqrt{2} }\log ||(z_1 +1, z_2 )||\, |<A \
 \ \mrm{  in  } \ \ \B ((-1,0), R) \ \nonumber \\
 &|g(z) - \frac{1}{\sqrt{2} }\log ||(z_1 -1, z_2 )||\,  |<A \ \
 \mrm{ in }\ \  \B ((1,0), R). \nonumber \end{eqnarray}

 Let $0<r<\frac{R}{16}$ be fixed and let $z_2=w$ be fixed with $0<|w| <r.$

 Consider the restriction $\hat g$: $\hat g _w (z)= \hat g (z,w)$. If $|z-1|\leq r$ or $|z+1|\leq r$ then $||(z,w)||< 2$ so $\hat g (z,w)\leq 0$ on
 $\{ |z-1|\leq r\}$ and $\{ |z+1|\leq r\}.$

 If  $-\infty\not\equiv\hat g _w \in \mc L _t  (\C  )$ one concludes that the total mass of $\Delta \hat g _w$ does not exceed $t .$ By
 symmerty one can assume that
\begin{equation}
\int _{\B (1,R)} \Delta \hat g _w \leq t/2 . \label{eq:4}
\end{equation}
(Otherwise consider $\B (-1,R)$ in place of $\B (1,R)$.)

If $|z-1|\leq |w|$ we then
have
\begin{equation}
\hat g _w (z) = \hat g(z,w)\leq \frac{1}{\sqrt{2} }\log ||(z-1, w)||+A  \\
 \leq \frac{1}{\sqrt{2} }\log |w| +A+1.  \label{eq:5}
\end{equation}

 Let  $z$ be any point on   $\{|w|<|z-1|\leq r\}$.
 Denote by $\B _1$ the disk $\B (z, 2r)$ and by $\B _2$ the disk
 $\B (1, r)$.
 Then $\B _2 \subset \B _1 \subset \B (1,R). $ If $J(K)$ denotes the
 average value of $\hat g _w$ over a set $K\subset \C$ then
\begin{equation}
\hat g _w (z)\leq J(\B _1 )\leq \frac{1}{4} J(\B _2).\nonumber
\end{equation}
Since $J(\B _2 )$ is dominated by the average of $\hat g _w$ over
the boundary of $\B _2$ one obtains from Riesz representation formula, using that  $\hat g _w \leq 0$   and
(\ref{eq:4}),
(\ref{eq:5}) :
\begin{eqnarray}
J(\B _2 ) &\leq \max _{\B (1, |w| ) } \hat g _w - \int _{\{|w|<|x-1|<r\}}
 \log|x-1|\Delta \hat g _w \nonumber \\
& \leq
\frac{1}{\sqrt{2} } \log |w| +A+1 - \frac{t}{2} \log |w|  \nonumber \\
& \leq ( \frac{1}{\sqrt{2}}-\frac{t}{2} )\log |w| +A+1 .
\nonumber
\end{eqnarray}

Therefore $$ \hat g (z,w) \leq  \frac{\sqrt{2}-t}{8} \log |w| +A +1$$ for
$||(z-1,w)|| < r .$ Since the Monge-Amp\`ere operator
cannot be defined for $v(z,w)=\log |w|$ it follows that the same
goes for the function $\hat g$. \fin

\begin{remark} The above example relies on a geometrical effect
which is also responsible for nonexistence of solutions to the \MA
equations in $\C\emph{P} ^n $ where we have on the right hand side
a generic combination of Dirac measures (cf. [Co]).
\end{remark}

\begin{remark} It follows from [S] that $\hat g (z,w) - \frac{\sqrt{2}-t}{8}log |w|$ is plurisubharmonic on  $\C^2.$
For an elementary proof, see [Ce4].
\end{remark}

\noindent Urban CEGRELL \\
Department of Mathematics \\
University of Umea \\
S-90187 Umea, Sweden \\

\noindent S\l awomir KO\L ODZIEJ \\
Jagiellonian University \\
Institute of Mathematics \\
\L ojasiewicza 6, 30-348 Krak\'ow, Poland \\

\noindent Ahmed ZERIAHI \\
Institut de Math\'ematiques de Toulouse \\
UPS, 118 Route de Narbonne \\
31062 Toulouse cedex, France.


\begin{thebibliography}{99}

 \bibitem [ACCP] {ACCP} {P. {\AA}hag, U. Cegrell, R. Czyz, R. and  H.
 H. Pham,}
 {\it Monge-Amp\`ere measures on pluripolar sets.} { J. Math. Pures Appl.}{ 92 (2009), 613--627.}


 \bibitem [BGZ] {BGZ} { S.Benelkourchi, V. Guedj, A. Zeriahi,} {\it Plurisubharmonic functions with weak singularities,}
Acta Universitatis Upsaliensis, Proceedings of the conference in honor of
C.Kiselman (Kiselmanfest, Uppsala, May 2006), 2009.

 \bibitem [BT1]{BT1} {E. Bedford and  B.A. Taylor,} { \it
 The Dirichlet problem for the complex Monge-Amp\`ere operator.}{ Invent. Math.}{ 37
 (1976), 1-44.}

 \bibitem [BT2]{BT2} {E. Bedford and  B.A. Taylor,}
 {\it  A new capacity for
 plurisubharmonic functions}, {Acta Math.}, {149 (1982)}, {1-40 }.

 \bibitem [BT3]{BT3} {E. Bedford and  B.A. Taylor,}
 {\it Fine topology, \v Silov boundary, and  $(dd^c)^n$},  J. Funct. Anal.  {72  (1987)},  no. 2, {225--251}.

\bibitem [Ce1] {Ce1} {U. Cegrell}, {\it Pluricomplex energy},
{Acta Math.}, { 180 (1998)}, {187-217}.

\bibitem [Ce2] {Ce2} {U. Cegrell}, {\it The general definition of
the complex Monge-Amp\`ere operator}, { Ann. Inst. Fourier}, { 51
(2004)}, {159-179}.


\bibitem [Ce3] {Ce3} { U. Cegrell}, {\it A general  Dirichlet problem
for the Complex Monge-Amp\`ere operator.} Ann. Pol. Math. 94.2 (2008), 131-147.


\bibitem [Ce4] {Ce4} { U. Cegrell}, {\it  Removable singularities for plurisubharmonic functions and related problems.}
       Proc. London Math. Soc. 36 (1978), 310-326.


\bibitem [CH] {CH} { U. Cegrell, L. Hed,} {\it Subextension and
approximation of negative plurisubharmonic functions.}  Michigan Math. J. Vol 56:3 (2008), 593-601.


\bibitem [CKZ] {CKZ} { U. Cegrell, S. Ko\l odziej and A. Zeriahi,}
{\it Subextension of plurisubharmonic functions with weak
singularities.} {Math. Z.} {250 (2005), 7-22}.

\bibitem [CZ] {CZ} { U. Cegrell and A. Zeriahi,}
{\it Subextension of plurisubharmonic functions with bounded Monge-Amp\`ere mass.} {C. R. Acad. Sci. Paris 336 (2003), no. 4, 305-308}.


\bibitem [Co] {Co} {D. Coman,} {\it Certain classes of pluricomplex Green functions on
$\C ^n$.}{ Math. Z.}  {235 (2000), 111-122}.

\bibitem [CGZ] {CGZ} {D. Coman, V. Guedj, A. Zeriahi,} {\it  Domains of definition of Monge-Amp\`ere operators on compact
K\"ahler manifolds}{ Math. Z.}  {259 (2008), 393-418}.


\bibitem [GZ1] {GZ1} {V. Guedj, A. Zeriahi,} {\it Intrinsic capacities on compact K\"ahler manifolds} {J. Geom. Anal. {\bf 15} (2005), no. 4, 607-639.}

\bibitem [GZ2] {GZ2} {V. Guedj, A. Zeriahi,} {\it The weighted Monge-Amp\`ere energy of
quasiplurisubharmonic functions.} {J. Funct. Anal.} {250 (2007),
442-482}.


 \bibitem [P] {P} {H.H. Pham}, {\it Pluripolar sets and the subextension in Cegrell's classes.} Complex Var. Elliptic Equ. 53(2008), 675-684.

 \bibitem [S] {S} {Y.T. Siu}, {\it Analyticity of sets associated to Lelong numbers and the extension of closed positive currents.}
Invent. Math. 27 (1974), 53-156.


\end{thebibliography}
 \end{document}